\documentclass[12pt]{article}
\pdfoutput=1
\usepackage{amssymb}
\usepackage{amsmath}
\usepackage{amsthm}
\usepackage{color}
\usepackage{comment}
\usepackage{geometry}		
\usepackage{graphicx}

\let\originalleft\left
\let\originalright\right
\renewcommand{\left}{\mathopen{}\mathclose\bgroup\originalleft}
\renewcommand{\right}{\aftergroup\egroup\originalright}

\usepackage{caption}
\usepackage{subcaption}

\usepackage{hyperref}
\hypersetup{
    colorlinks=true,
    linkcolor=black,
    urlcolor=blue,
    citecolor=black
}

\begin{document}

\newcommand\cO{\mathcal{O}}
\newcommand{\ee}{\varepsilon}

\newcommand{\removableFootnote}[1]{\footnote{#1}}

\newtheorem{theorem}{Theorem}[section]
\newtheorem{corollary}[theorem]{Corollary}
\newtheorem{lemma}[theorem]{Lemma}
\newtheorem{proposition}[theorem]{Proposition}

\theoremstyle{definition}
\newtheorem{definition}{Definition}[section]
\newtheorem{example}[definition]{Example}

\theoremstyle{remark}
\newtheorem{remark}{Remark}[section]



\newcommand{\orcid}[1]{\href{https://orcid.org/#1}{\includegraphics[width=8pt]{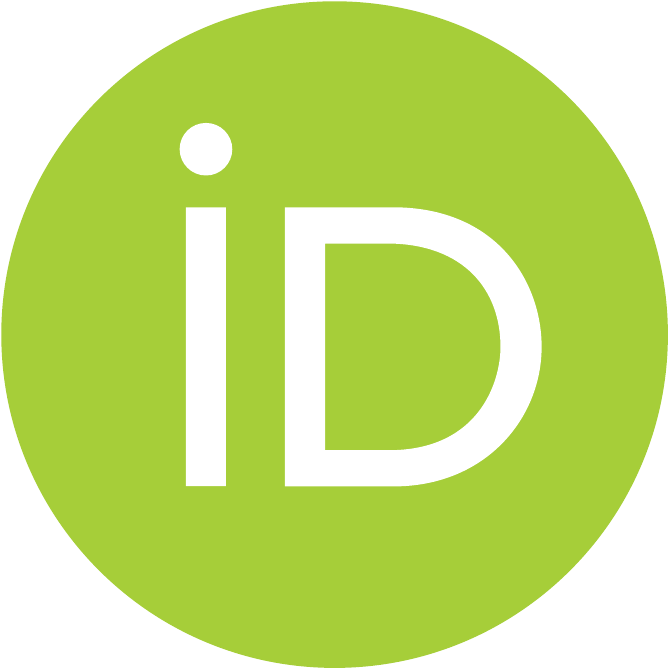}}}

\newcommand\nnfootnote[1]{%
  \begin{NoHyper}
  \renewcommand\thefootnote{}\footnote{#1}%
  \addtocounter{footnote}{-1}%
  \end{NoHyper}
}

\title{
Numerical bifurcation analysis of improved denatured Morris-Lecar neuron model
}
\author{
H.O.~Fatoyinbo\textsuperscript{a*}\orcid{0000-0002-6036-2957}, S.S.~Muni\textsuperscript{a}\orcid{0000-0001-9545-8345}, I.~Ghosh\textsuperscript{a}, I.O.~Sarumi\textsuperscript{b}, A. Abidemi\textsuperscript{c}\orcid{0000-0003-1960-0658}\\\\
\textsuperscript{a}\textit{School of Mathematical and Computational Sciences} \\\textit{Massey University, Palmerston North, New Zealand}\\
\textsuperscript{b}\textit{Department of Mathematics and Statistics} \\
\textit{King Fahd University of Petroleum and Minerals, Dhahran, Saudi Arabia}\\\textsuperscript{c}\textit{Department of Mathematical Sciences} \\
\textit{Federal University of Technology, Akure, Nigeria}
 
 }
 
 \date{}

\maketitle


\begin{abstract}
It is well-known that the electrical activities of neurons are induced by a wide variety of external factors. This work considers the effect of electromagnetic induction on improved denatured Morris-Lecar neuron model. The dependence of dynamical behaviour of the original denatured Morris-Lecar model on parameters is addressed through numerical bifurcation analysis. This allows us to explore the changes in dynamics of the model qualitatively as parameters are varied. Then we investigate the effects of external periodic current and electromagnetic flux on the dynamical properties of the improved denatured Morris-Lecar neuron model. Different types of dynamical behaviour, ranging from regular periodic spiking to complex bursting, are found when the multiple parameters are varied simultaneously. The improved model could be applied to research where simple models are required for physiological and pathophysiological responses in neurons.
\end{abstract}

\nnfootnote{E-mail addresses: H. O. Fatoyinbo (\url{h.fatoyinbo@massey.ac.nz}), S.S. Muni (\url{s.muni@massey.ac.nz}), I. Ghosh(\url{i.ghosh@massey.ac.nz}), I. O. Sarumi(\url{isarumi@kfupm.edu.sa}), A. Abidemi (\url{aabidemi@futa.edu.ng})}

\section{Introduction}\label{sec:intro}
Nonlinear dynamical behaviours of neurons have been widely studied to understand the mechanisms associated with neuronal diseases such as epilepsy \cite{ayala73} and Parkinson’s disease \cite{Levy2000}. They are frequently modelled by a nonlinear system of differential equations. In 1952, Hodgkin and Huxley 
(H-H) constructed a neuronal model from experimental studies to describe the conduction and propagation of electrical impulses in the squid giant axon \cite{Hodgkin1948TheAxon}. Following the H-H model, FitzHugh and Nagumo \cite{Fitzhugh}, \cite{Nagumo1962} constructed a two-variable model which is a simplified version of the H-H model to describe the different phenomena in neuron dynamics. Other well-known models are the Morris and Lecar (ML) model \cite{Morris}, the Hindmarsh–Rose model \cite{Indmarsh11984}, and the Izhikevich neuron model \cite{Izhikevich2007DynamicalBursting}.

Bifurcation analysis has been employed in many previous studies to uncover transitions of firing patterns in neurons and other biological cells \cite{Govaerts2005TheApproach}, \cite{Duan2006Codimension-twoModel}, \cite{hammed2}. It is very helpful to identify important parameters that can cause significant change in the dynamical behaviour qualitatively. Firing activities of neurons are classified into different classes of excitability depending on the transitions from a rest state to a firing (oscillatory) state. The authors in \cite{Rinzel1999AnalysisNetwork} classified neurons into two classes of excitability based on the onsets of firing. In class I excitability, the periodic oscillations occur through a saddle-node on an invariant circle (SNIC) bifurcation. In contrast, for class II excitability the periodic oscillations emerge through a Hopf bifurcation. The mode of transitions between the two classes of excitability is characterised by bifurcations. On this basis, the codimension-two bifurcations associated with a change between classes of excitability are identified in previous studies of neurons and other excitable cells \cite{Tsumoto2006BifurcationsModel}, \cite{Duan2008Two-parameterModel}, \cite{hammed}, \cite{WU2021}.

Multiple researchers have studied various factors that influence the firing activities in neurons and other excitable cells \cite{Zhao2017TransitionsAutapse},\cite{Liu2014BifurcationModelb},\cite{Govaerts2005TheApproach},\cite{Gonzalez-Miranda2014}. The effect of direct current stimulus \cite{Hetongwang2011}, temperature \cite{Xing2020}, pressure \cite{hammed}, time-delay \cite{Jia2018_2} and  environmental noise \cite{Gong2011} have been investigated. Magnetic induction is another factor that affects the neuron membrane potential significantly, which has been observed in physiological experiments \cite{Barry2016},\cite{eteme2019}. According to Maxwell's electromagnetic induction theorem, variation of neuronal membrane potential results in the generation of a magnetic field in the neuron and the induction current can influence the firing activities of the neuron. Li {\em et al.} \cite{Li2015} constructed a mathematical model to study the variation of firing patterns of neurons under the influence of electromagnetic radiation. Lv and Ma \cite{Lv2016} used the magnetic flux to describe the effects of electromagnetic radiation on the electrical activities of neurons by coupling the magnetic flux with the membrane potential. Many recent papers have studied the effects of electromagnetic induction on neurons and collective behaviours in a network of neurons \cite{Yang2021}, \cite{muni2022}, \cite{Mondal2019BifurcationModel}, \cite{Rajagopal2021}. Diverse firing modes and complicated dynamical behaviour including bursting are observed in neuron models under electromagnetic radiation. Also, spatiotemporal dynamics including chaotic bursting, pattern formation and chimera state are found in a network of neurons \cite{KARTHIKEYAN2021}.

In many studies, high-dimensional conductance-based models have been used to explore electrophysiological activities in neurons. However, due to the nonlinear behaviour of the models, they are computationally expensive; hence low-dimensional models have been instrumental to study electrical activities in neurons because the reduced models capture the excitable properties in higher-dimensional models. A two-variable ML model is a reduced model that replicates the dynamical behaviours observed in the 8-dimensional H-H model. In \cite{Schaeffer}, a simplified version of the ML model, known as the denatured ML neuron model, is introduced. To the best of our knowledge, this is the first work to give a more comprehensive analysis of the denatured ML model

In this paper, we will use the denatured ML model to investigate the dynamics of neurons under the influence of electromagnetic induction and external periodic current. The paper is structured as follows: In Section~\ref{sec:model}, we describe the model and the parameters. In Section~\ref{sec:DML_original}, we analyse the original denatured ML model via phase plane and bifurcation analyses to explore the firing patterns of neurons. We also explore the effects of electromagnetic flux and external periodic current on firing patterns of the improved ML model in Section~\ref{sec:DML_improved}. The summary and conclusion of results are given in Section~\ref{sec:conclusion}.

\section{Model Description}\label{sec:model}
The denatured Morris-Lecar neuron model is a simplified version of the two variable Morris-Lecar neuron model \cite{Morris}. It consists of two-coupled ODEs described in \cite{Schaeffer} as
\begin{subequations}
\begin{align}
 \dot{x} & \!\begin{aligned}[t]&=x^2(1-x)-y+I, \end{aligned}\\[0.75ex]
 \dot{y} & \!\begin{aligned}[t]&=Ae^{\alpha x}-\gamma y, \end{aligned}
\end{align}
\label{eq:DML}
\end{subequations}
where $x$ is the membrane action potential and $y$ is the recovery variable. The system parameter $I$ is the external current, $A$, $\alpha$ and $\gamma$ are positive constants. In this paper we improve model \eqref{eq:DML} by adding electromagnetic induction and external periodic current. The dynamical equations of the improved version of \eqref{eq:DML} are described as follows
\begin{subequations}
\begin{align}
 \dot{x} & \!\begin{aligned}[t]&=x^2(1-x)-y+I+k\rho(\phi)x, \end{aligned}\\[0.75ex]
 \dot{y} & \!\begin{aligned}[t]&=Ae^{\alpha x}-\gamma y, \end{aligned}\\[0.75ex]
\dot{\phi}&\!\begin{aligned}[t]&=k_{1}x-k_{2}\phi+\phi_{\rm ext},  \end{aligned}
\end{align}
\label{eq:DML_em}
\end{subequations}
 where $k$ is the feedback gain, the function $\rho(\phi)=\alpha_{1}+3\beta\phi^{2}$ is the electromagnetic effect on the neuron and $\phi$ is the magnetic flux across the cell membrane. $\alpha_{1}$ and $\beta$ are the memory conductances, and $\phi_{\rm ext}$ is the external magnetic flux. The external periodic current $I=I_{0}\sin(\omega t)$, where $I_{0}$ is the amplitude of the current and $\omega$ is the angular frequency. Unless otherwise stated, the parameter values except for $I_{0}$ and $\gamma$ are fixed as:
 $A=0.0041$, $\alpha=5.276$, $k=0.003$, $k_{1}=0.19$, $k_{2}=0.5$, $\alpha_{1}=0.1$, $\beta=0.02$, $\omega=0.01$, and $\phi_{\rm ext}=0.2$.

\section{Analysis of the original denatured Morris-Lecar model}\label{sec:DML_original}
\begin{figure}
    \centering
    \includegraphics[scale=0.45]{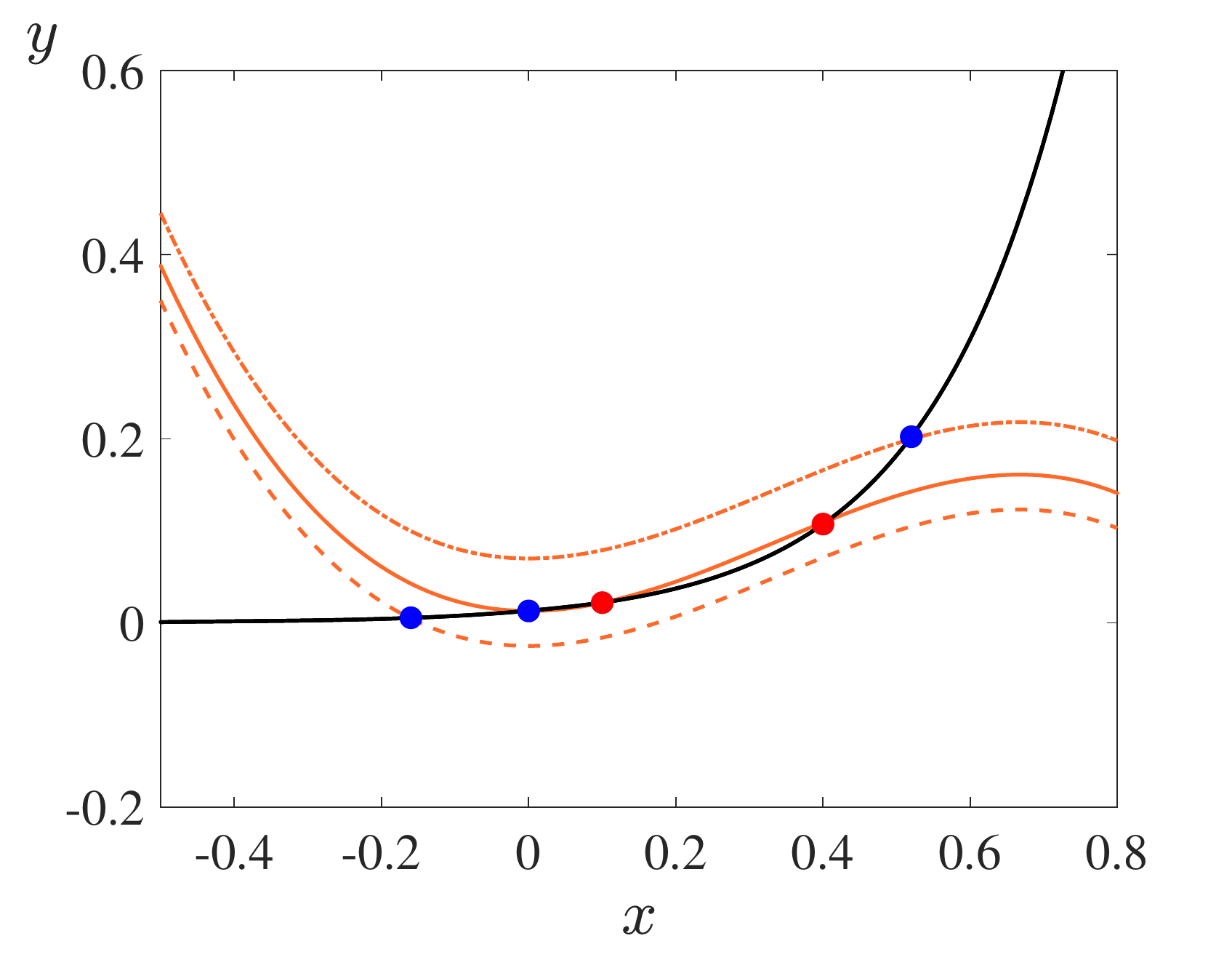}
    \caption{The $(x,y)$ phase plane of \eqref{eq:DML}. The $y$-nullcline is plotted black and the $x$-nullclines are plotted orange for three different values of $I$. The dashed, solid and dashed-dot curves are $x$-nullclines for $I=-0.025$, $0.013$, and $0.07$, respectively. The blue [red] filled circle correspond to stable and unstable equilibria.}
    \label{fig:nullcline}
\end{figure}

We begin our analysis by investigating the dynamical behaviour of the original denatured ML system \eqref{eq:DML} via the interaction of the variables in the phase plane. The $x$- and $y$-nullclines are set of points in the phase plane such that $\dot{x}=0$ and $\dot{y}=0$, that is $x$ and $y$ are constant, and the intersection of the nullclines defines the equilibrium point. The $y$-nullcline is shown in black and the $x$-nullclines are shown in orange for various values of $I$ in Fig.~\ref{fig:nullcline}. The blue [red] filled circle corresponds to stable [unstable] equilibrium. For extremely low values of $I$ the two nullclines intersect at one point, corresponding to a unique stable equilibrium. For example, the dashed orange curve in Fig.~\ref{fig:nullcline} is when $I=-0.025$. For intermediate values of $I$ the nullclines intersect at three points, corresponding to three equilibria of the system. The case of $I=0.013$ is the solid orange curve in Fig.~\ref{fig:nullcline}, there is one stable and two unstable equilibria. For large values of $I$ the two nullclines intersect at one point, this corresponds to a unique equilibrium point. The point is either stable or unstable depending on the value of $I$. For $I=0.07$, the equilibrium point is stable, see the dash-dot orange curve in Fig.~\ref{fig:nullcline}. From the analysis above, it is observed that the number and stability of equilibria of the system change as the external current $I$ is varied, thus to further our understanding of how the model dynamics changes as parameters are varied we will consider numerical bifurcation analysis of system \eqref{eq:DML} in the next section.  

\subsection{One-parameter bifurcation analysis}\label{sec:codim-1}
In this section, we consider the dependence of the dynamical behaviour of \eqref{eq:DML} on parameters $I$ and $\gamma$ through numerical bifurcation analysis. The bifurcation diagrams were produced in MatCont, a numerical continuation software \cite{Govaerts2011MATCONTMatlab}. Fig.~\ref{fig:(I,x)} is a one-parameter bifurcation diagram of the model \eqref{eq:DML} for the bifurcation parameter $I$. The blue [red] curves correspond to stable [unstable] solutions, respectively. The points labelled ${\rm LP}_{1}$ and ${\rm LP}_{2}$ are saddle-node bifurcation points. The points labelled ${\rm HB}_{1}$ and ${\rm HB}_{2}$ are Hopf bifurcation points. For relatively low values of $I$ there exists a unique stable solution. With increasing $I$, the stable solution loses stability in a saddle-node bifurcation ${\rm LP}_{1}$ at $I=0.0153$. Upon further increase in the value of $I$, the unstable solution regains stability in a saddle-node bifurcation ${\rm LP}_{2}$ at $I=0.0109$. Between the two saddle-node bifurcations, ${\rm LP}_{1}$ and ${\rm LP}_{2}$, there exist three solutions: two stable and one unstable. Thus, in this parameter regime there is bistability (two stable solutions coexist) in the system. The stable solution produced at ${\rm LP}_{2}$ loses stability in a Hopf bifurcation ${\rm HB}_{1}$ at $I=0.0117$. Between saddle-node ${\rm LP}_{1}$ and Hopf bifurcation ${\rm HB}_{1}$ there are three solutions: one stable (lower branch) and two unstable (middle and upper branch). An example of this scenario is shown in Fig.~\ref{fig:nullcline} when $I=0.013$. Since the first Lyapunov coefficient of ${\rm HB}_{1}$ is negative ($a=-5.2312$) then ${\rm HB}_{1}$ is a supercritical bifurcation. As known from the theory of supercritical Hopf bifurcation \cite{KuznetsovY.A.1995ElementsTheory}, a stable periodic solution emanates from ${\rm HB}_{1}$.  As the value of $I$ increases the stable periodic solution terminates in a Hopf bifurcation ${\rm HB}_{2}$ at $I=0.0541$. The first Lyapunov coefficient of ${\rm HB}_{2}$ is negative ($a=-0.2643$), thus it is a supercritical Hopf bifurcation. Also, the unstable solution from ${\rm HB}_{1}$ regains stability at ${\rm HB}_{2}$. For relatively high values of $I$ there is a unique stable solution.
\begin{figure}[!htbp]
\centering
   \begin{subfigure}[b]{.45\linewidth}
    \centering
    \caption{}
    \includegraphics[width=.99\textwidth]{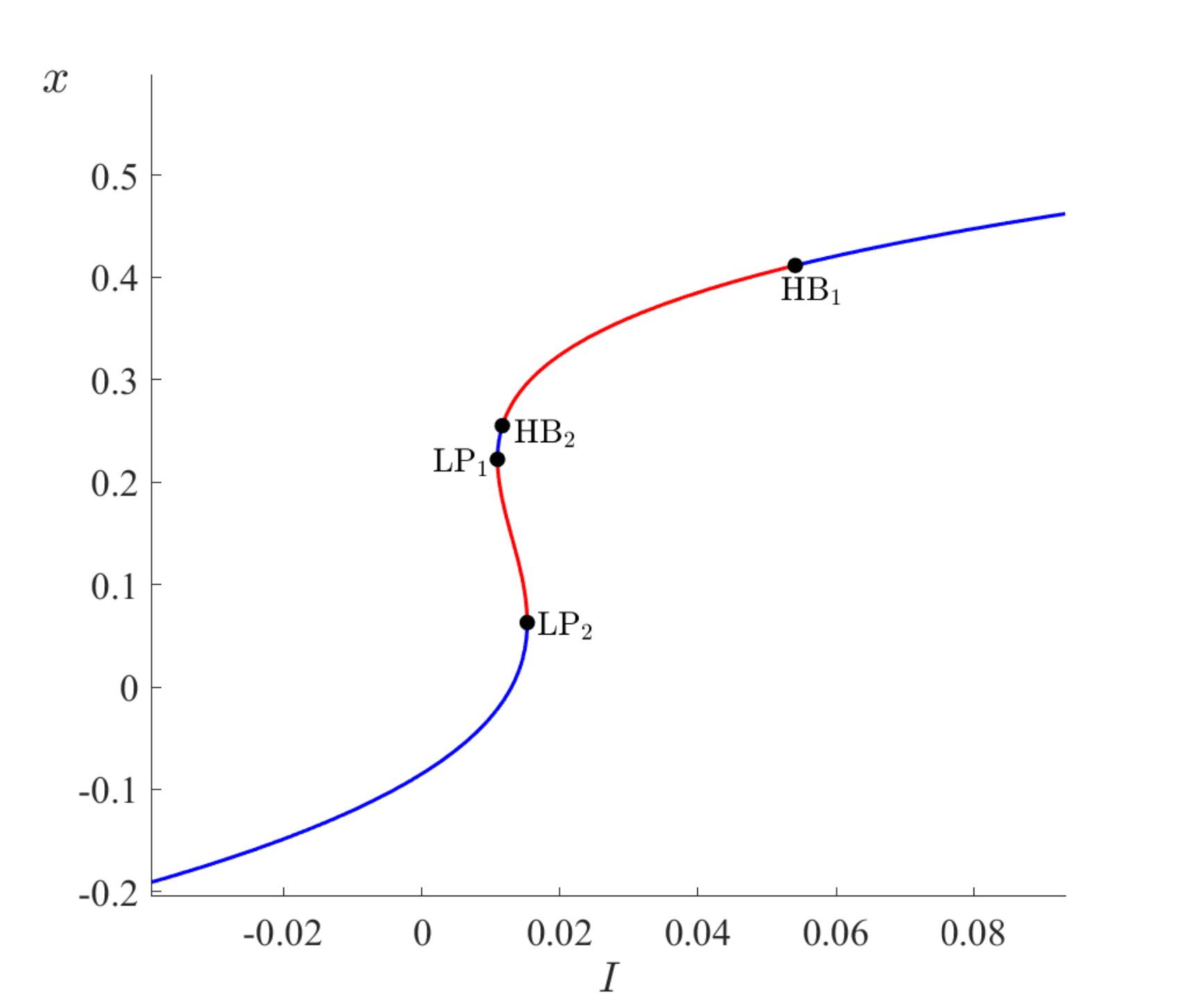}
     \label{fig:(I,x)}
  \end{subfigure}%
  \begin{subfigure}[b]{.45\linewidth}
    \centering
    \caption{}
    \includegraphics[width=.99\textwidth]{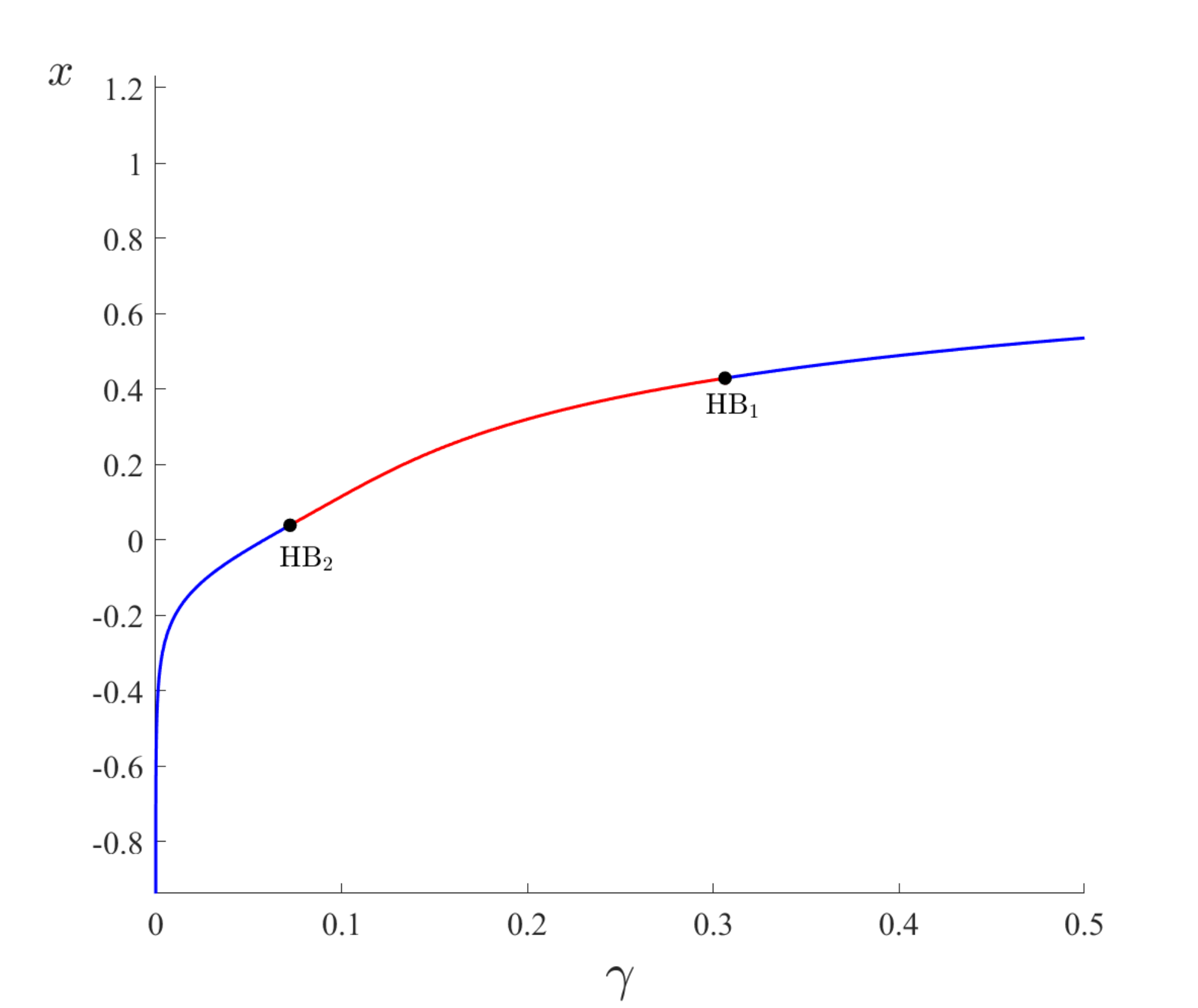}
     \label{fig:(gamma,x)}
  \end{subfigure}%
 \caption{Bifurcation diagrams of \eqref{eq:DML} with (a) $I$; (b) $\gamma$ as bifurcation parameters. Blue [red] curves correspond to stable [unstable] solution. HB: Hopf bifurcation; LP: saddle-node bifurcation
(of an equilibrium).}
 \label{fig:bifurcations_c}
\end{figure}

Next $\gamma$ is varied. A bifurcation diagram is shown in Fig.~\ref{fig:(gamma,x)}. We observed that for the parameter values of $\gamma$ considered, the system has a unique solution. Apart from the parameter regime between the two Hopf bifurcations ${\rm HB}_{1}$ and ${\rm HB}_{2}$ where the solution is unstable, it is stable everywhere else (to the left ${\rm HB}_{1}$ and the right of ${\rm HB}_{2}$). The first Lyapunov coefficients of ${\rm HB}_{1}$ and ${\rm HB}_{2}$ are negative hence they are supercritical Hopf bifurcations. The stable periodic solution that emanates from ${\rm HB}_{1}$ terminates at ${\rm HB}_{2}$. From the above analysis, system \eqref{eq:DML} models neuron activity with class II excitability.

\subsection{Two-parameter bifurcation analysis}\label{sec:codim-2}
\begin{figure}[!htbp]
    \centering
    \includegraphics[scale=0.4]{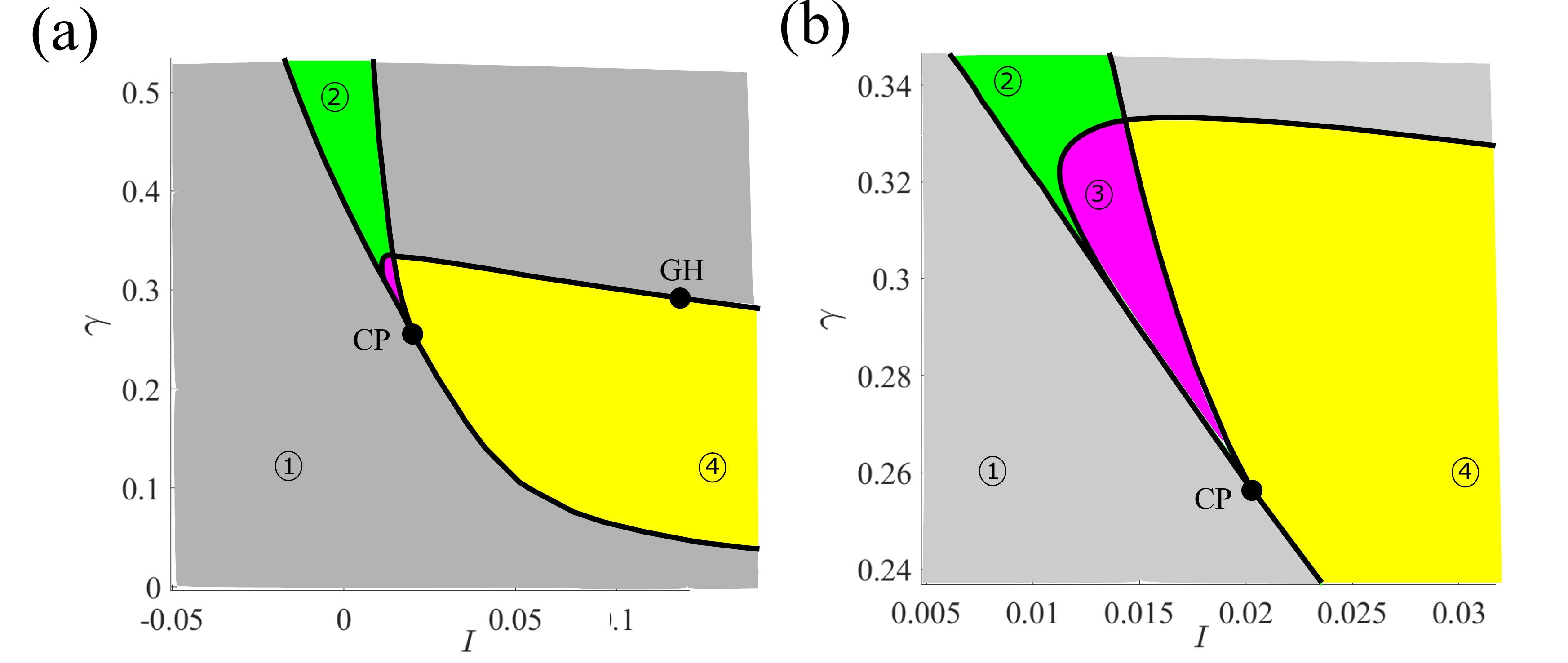}
    \caption{A two-parameter bifurcation diagram of model \eqref{eq:DML} in the $(I,\gamma)$-plane; (b) is the enlargement of Fig.~\ref{fig:Two-par}a. The black curves are the loci of the saddle-node and Hopf bifurcations in Fig.~\ref{fig:bifurcations_c}. The codimension-2 bifurcation points are CP: cusp bifurcation; GH: generalised Hopf bifucation.}
    \label{fig:Two-par}
\end{figure}

Here a two-parameter bifurcation analysis of model \eqref{eq:DML} by varying $I$ and $\gamma$ simultaneously is considered. Fig.~\ref{fig:Two-par} is a codimension-2 bifurcation diagram of \eqref{eq:DML} in $(I,\gamma)$-plane. The black curves correspond to the loci of the Hopf and saddle-node bifurcations in Fig.~\ref{fig:(I,x)}. The saddle-node loci collide and annihilate in a codimension-2 bifurcation known as cusp bifurcation CP. A generalised Hopf bifurcation GH appears along the Hopf locus, this is a codimesion-2 point where the Hopf bifurcation  changes from subcritical to supercritical. Fig.~\ref{fig:Two-par} is divided into four qualitatively different dynamical behaviour regions, with enlargement in Fig.~\ref{fig:Two-par}b. Each region is assigned a number and a colour. The remainder of this section describes the dynamical behaviour in each region.
In region \textcircled{\scriptsize{1}}, the system has one stable equilibrium and no periodic solutions (rest state). In region \textcircled{\scriptsize{2}}, the system has three equilibria: two stable and one unstable equilibria. This region is bounded by the loci of saddle-node bifurcations, thus there is no periodic solution. Apart from the loci of the saddle-node bifurcations previously observed there is now a Hopf bifurcation locus between which there exists a stable periodic solution in region \textcircled{\scriptsize{3}}. Finally, in region \textcircled{\scriptsize{4}} there is one unstable equilibrium and one stable periodic solution.

\section{Analysis of improved denatured Morris-Lecar model}\label{sec:DML_improved}
\begin{figure}[htbp]
\centering
\begin{subfigure}[b]{.26\textwidth}
  \centering
  \caption{}
  \includegraphics[width =\textwidth]{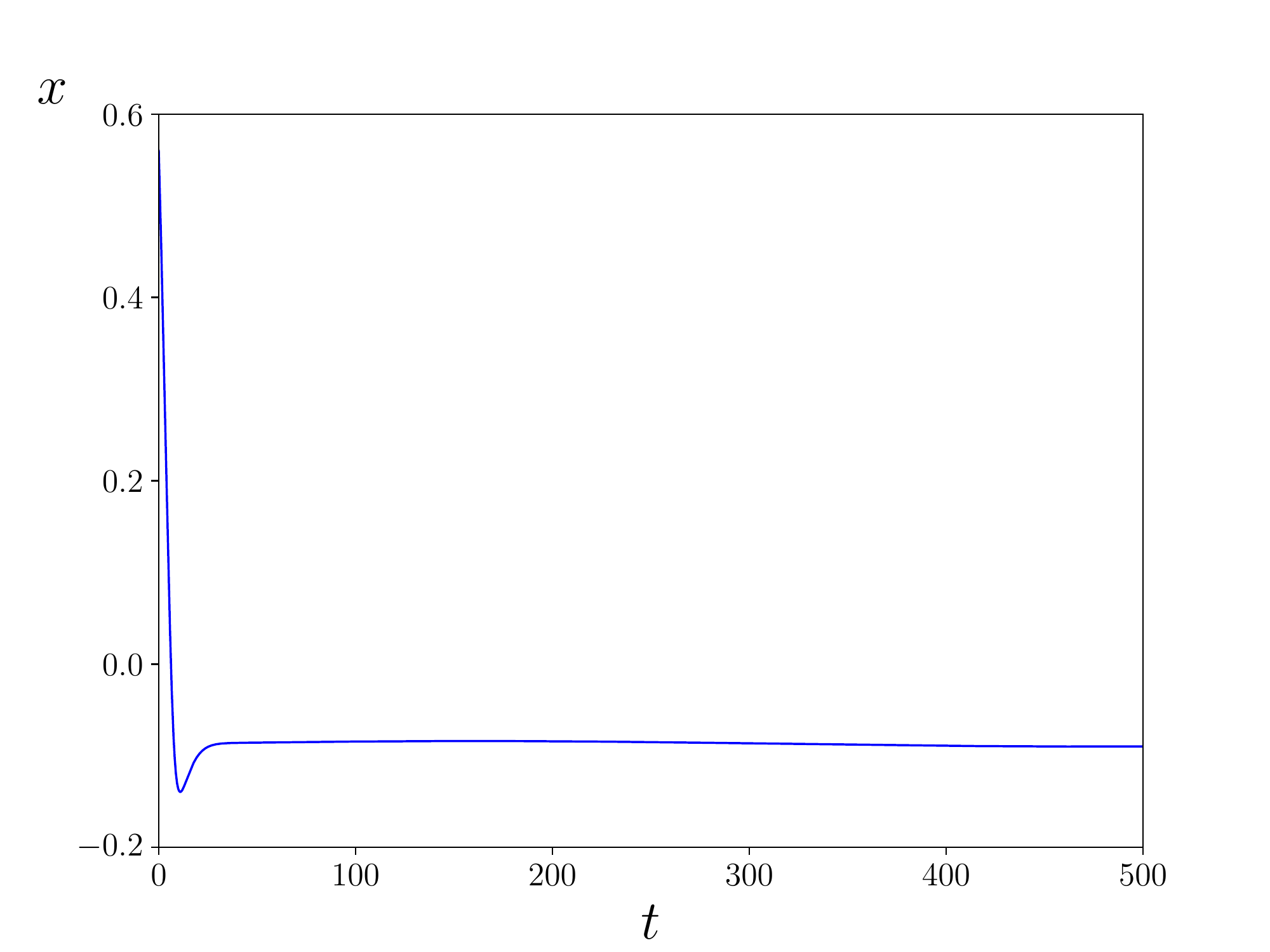}
  \label{fig:0.00072}
\end{subfigure}%
\begin{subfigure}[b]{.26\textwidth}
  \centering
  \caption{}
  \includegraphics[width =\textwidth]{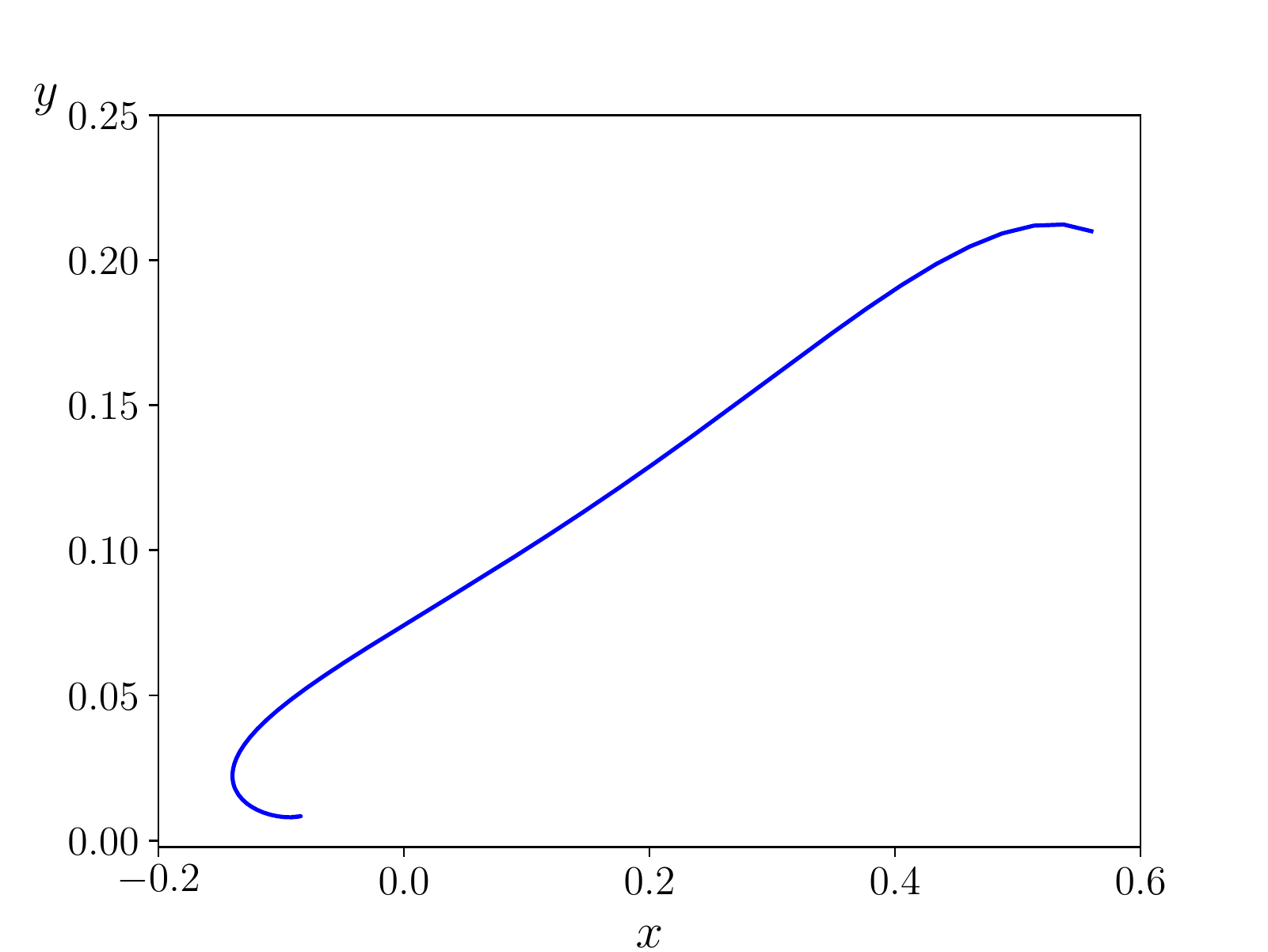}
  \label{fig:0.00072_pp}
\end{subfigure}%
\begin{subfigure}[b]{.26\linewidth}
    \centering
    \caption{}
    \includegraphics[width=\textwidth]{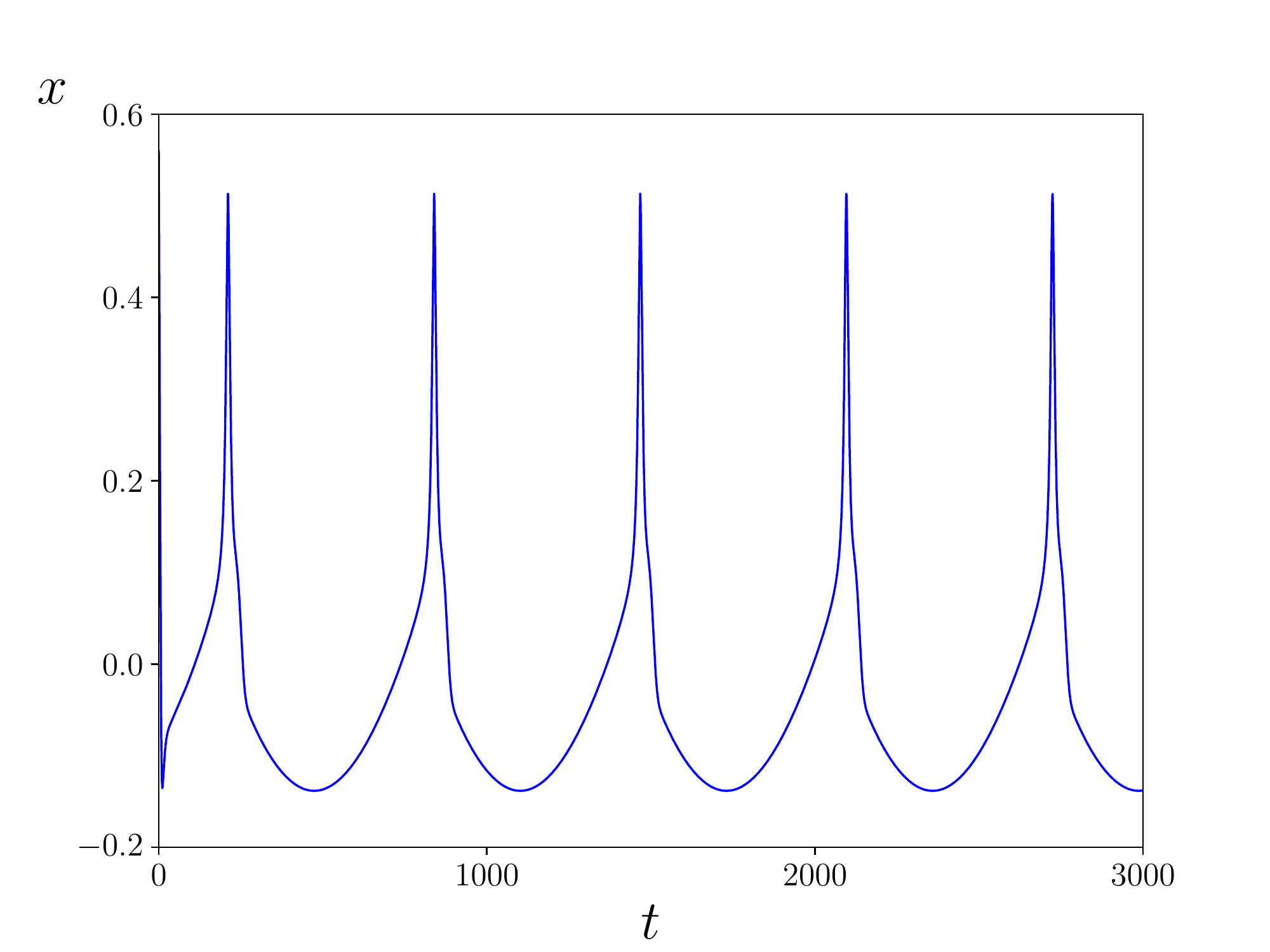}
     \label{fig:0.0155}
  \end{subfigure}%
  \begin{subfigure}[b]{.26\linewidth}
    \centering
    \caption{}
    \includegraphics[width=\textwidth]{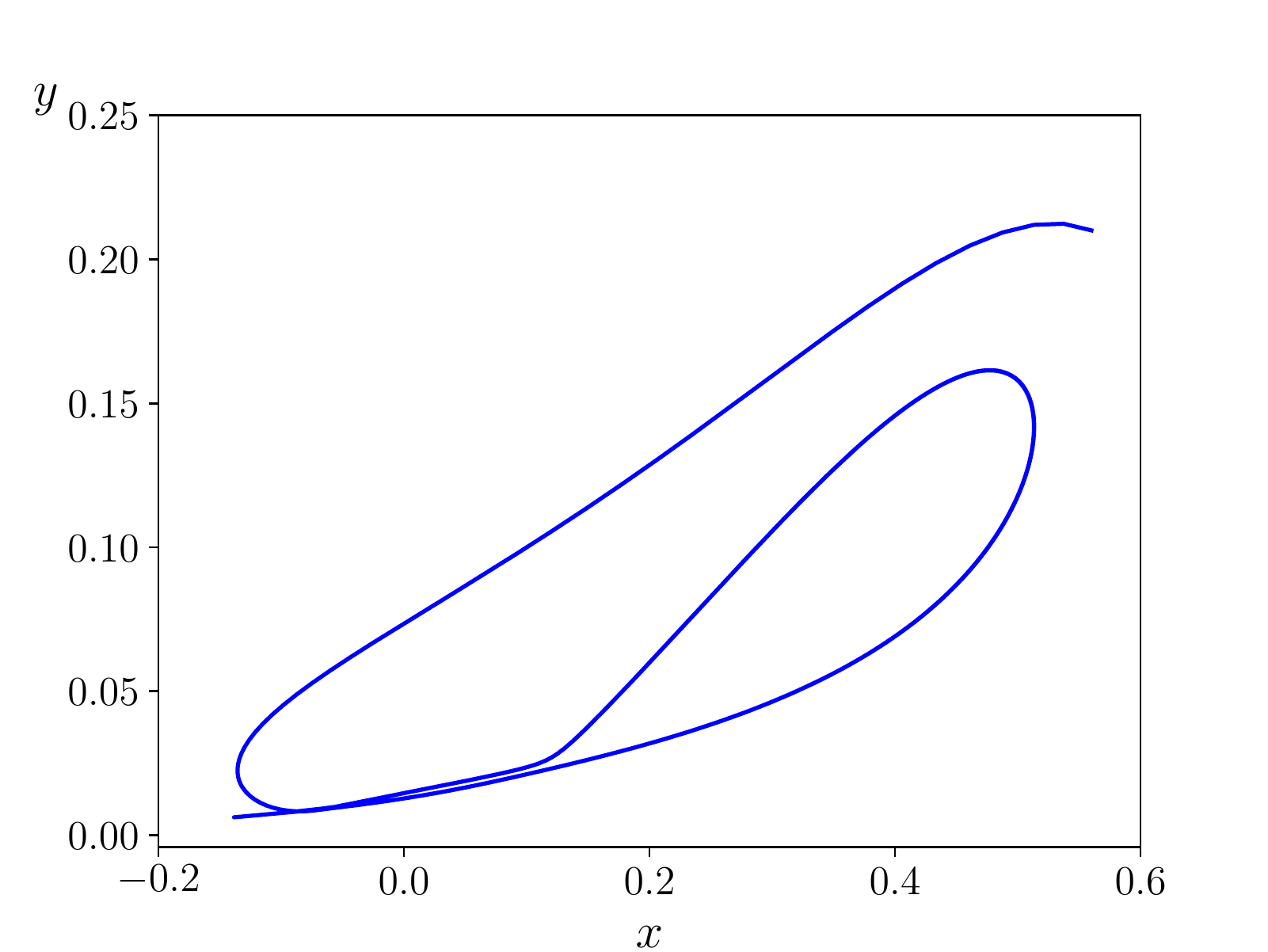}
     \label{fig:0.0155_pp}
  \end{subfigure}\\%
  \begin{subfigure}[b]{.26\textwidth}
  \centering
  \caption{}
 \includegraphics[width =\textwidth]{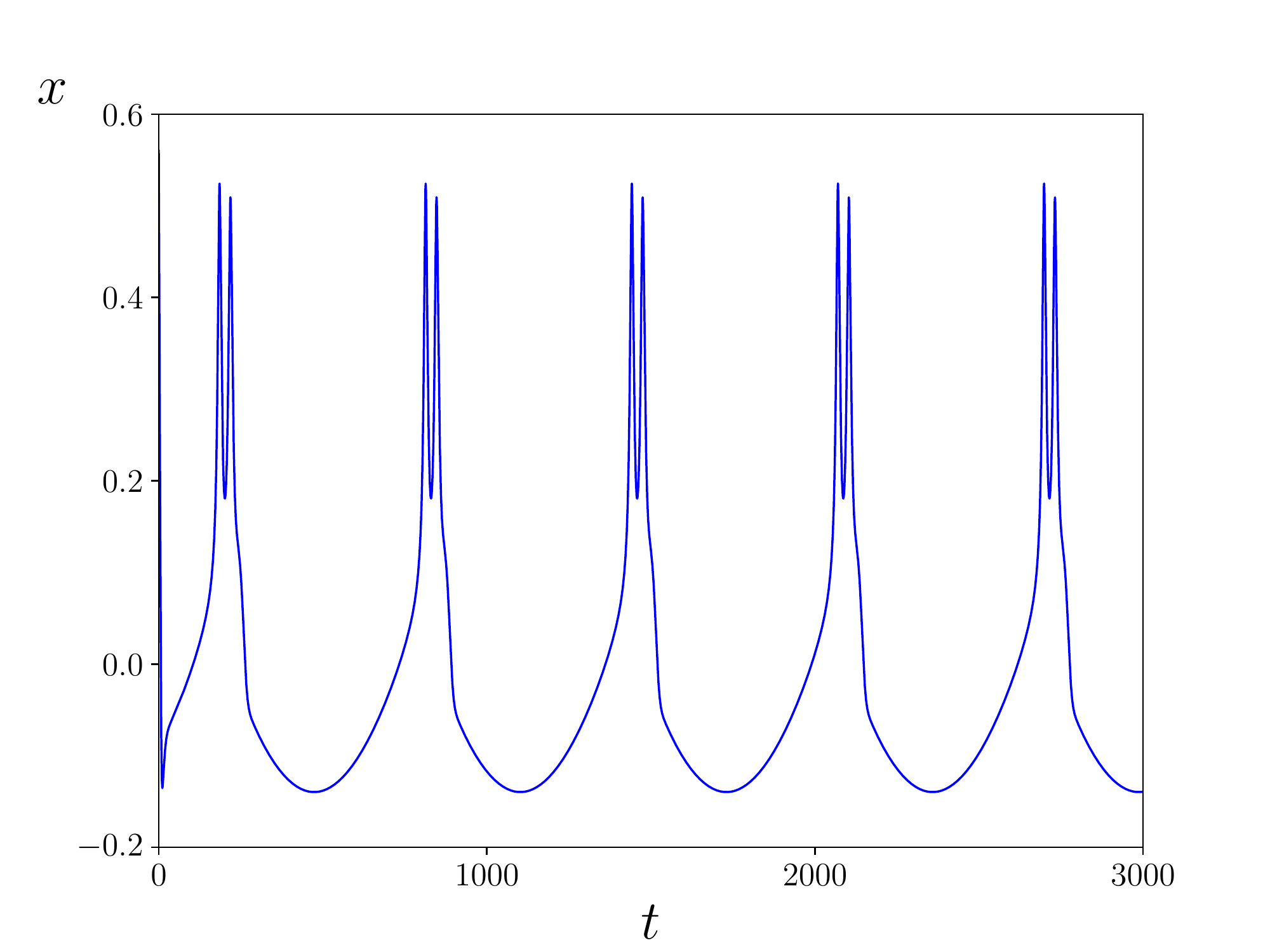}
  \label{fig:0.016}
\end{subfigure}%
\begin{subfigure}[b]{.26\textwidth}
  \centering
  \caption{}
 \includegraphics[width = \textwidth]{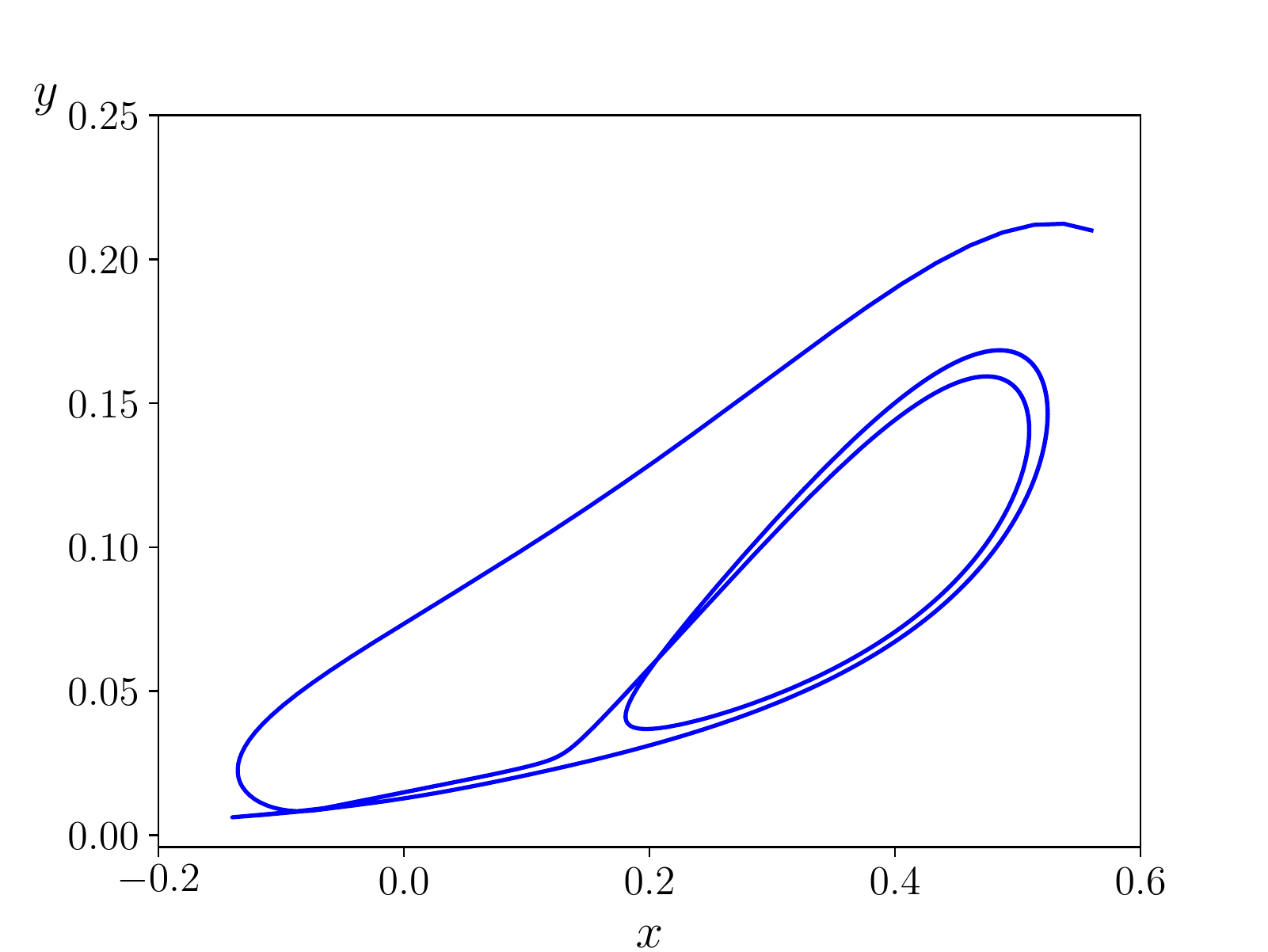}
  \label{fig:0.016_pp}
\end{subfigure}%
\begin{subfigure}[b]{.26\textwidth}
  \centering
  \caption{}
 \includegraphics[width = \textwidth]{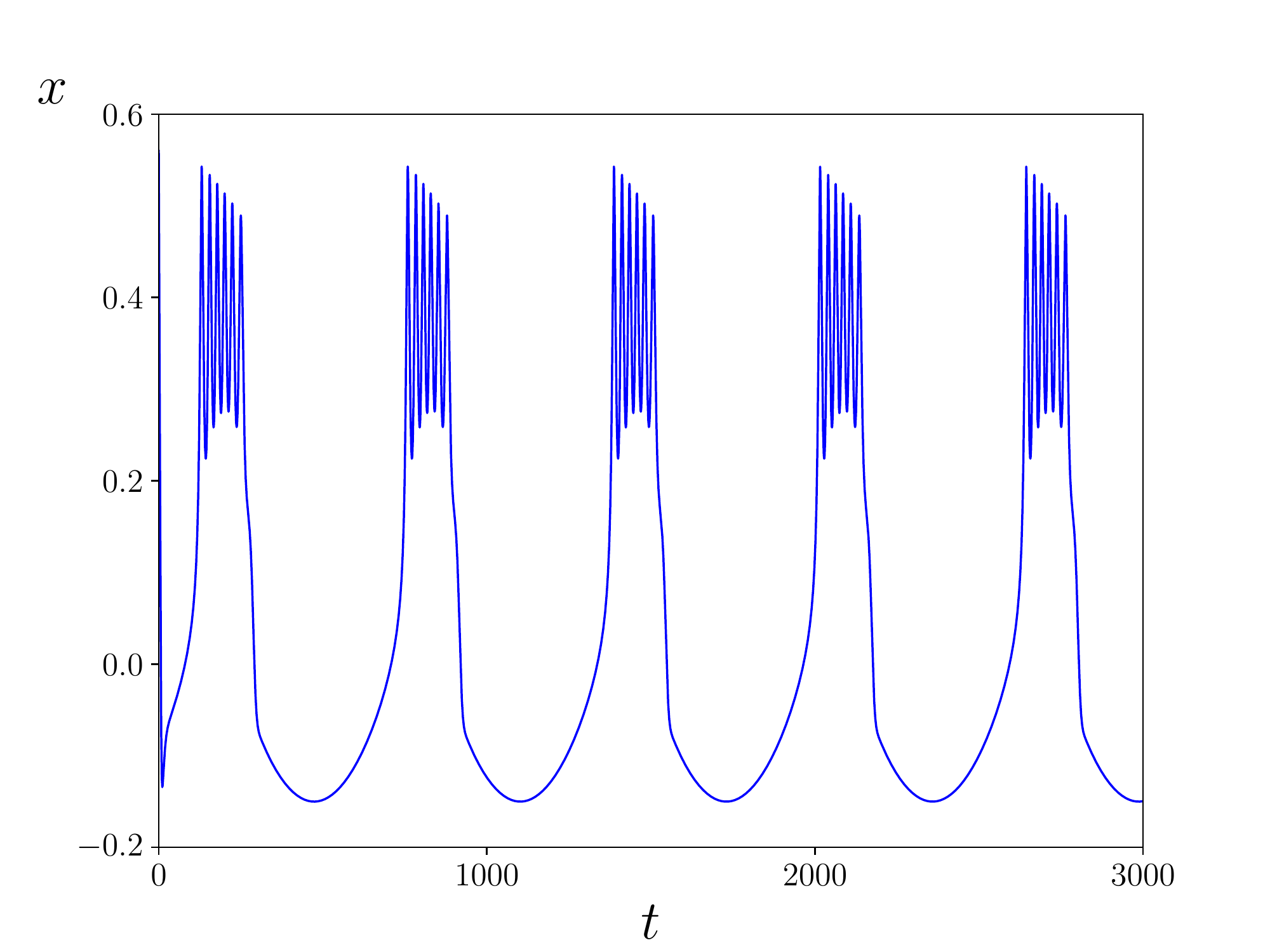}
  \label{fig:0.02}
\end{subfigure}%
\begin{subfigure}[b]{.26\textwidth}
  \centering
  \caption{}
 \includegraphics[width = \textwidth]{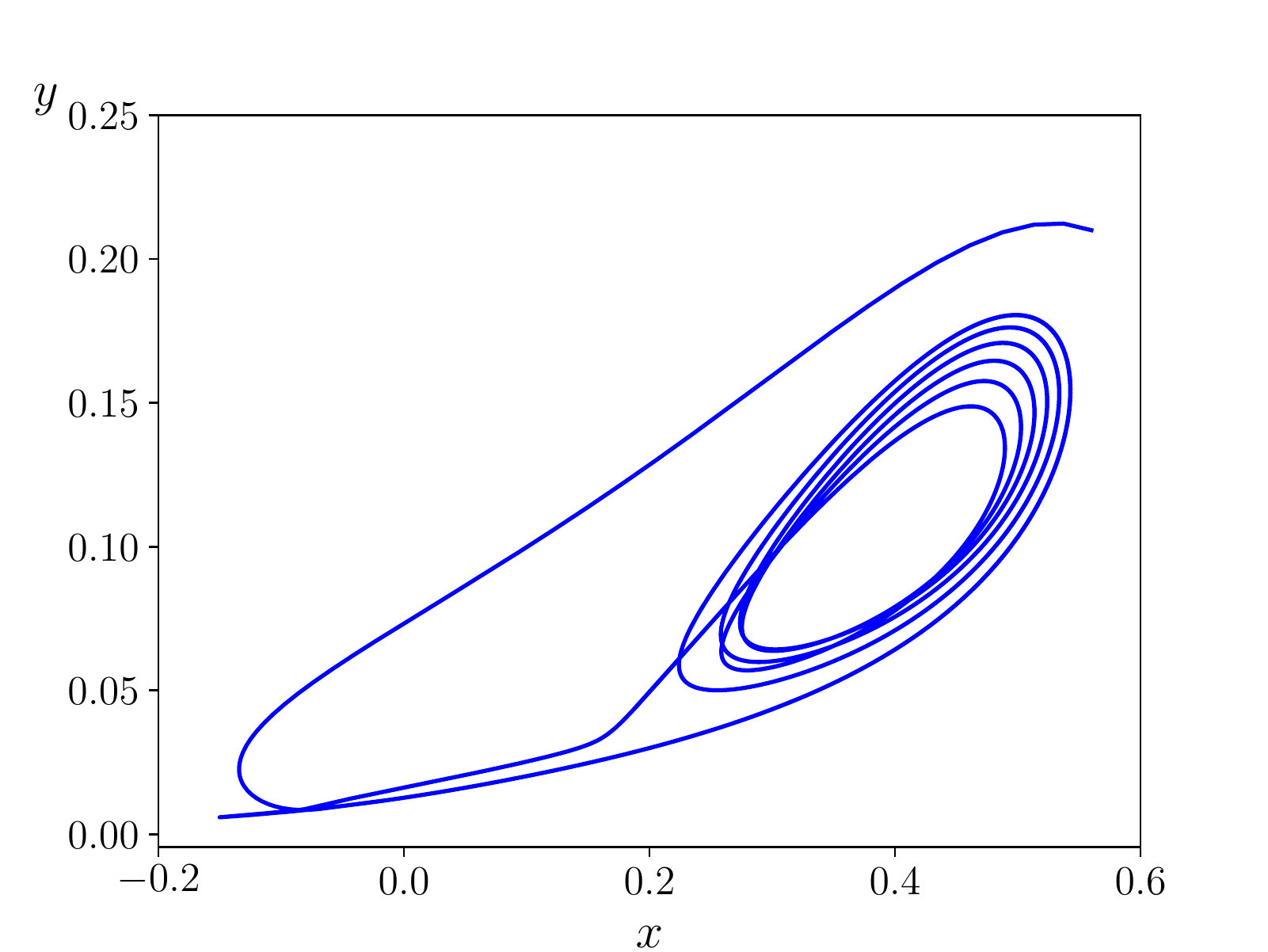}
  \label{fig:0.02_pp}
\end{subfigure}\\%
\begin{subfigure}[b]{.26\textwidth}
  \centering
  \caption{}
 \includegraphics[width = \textwidth]{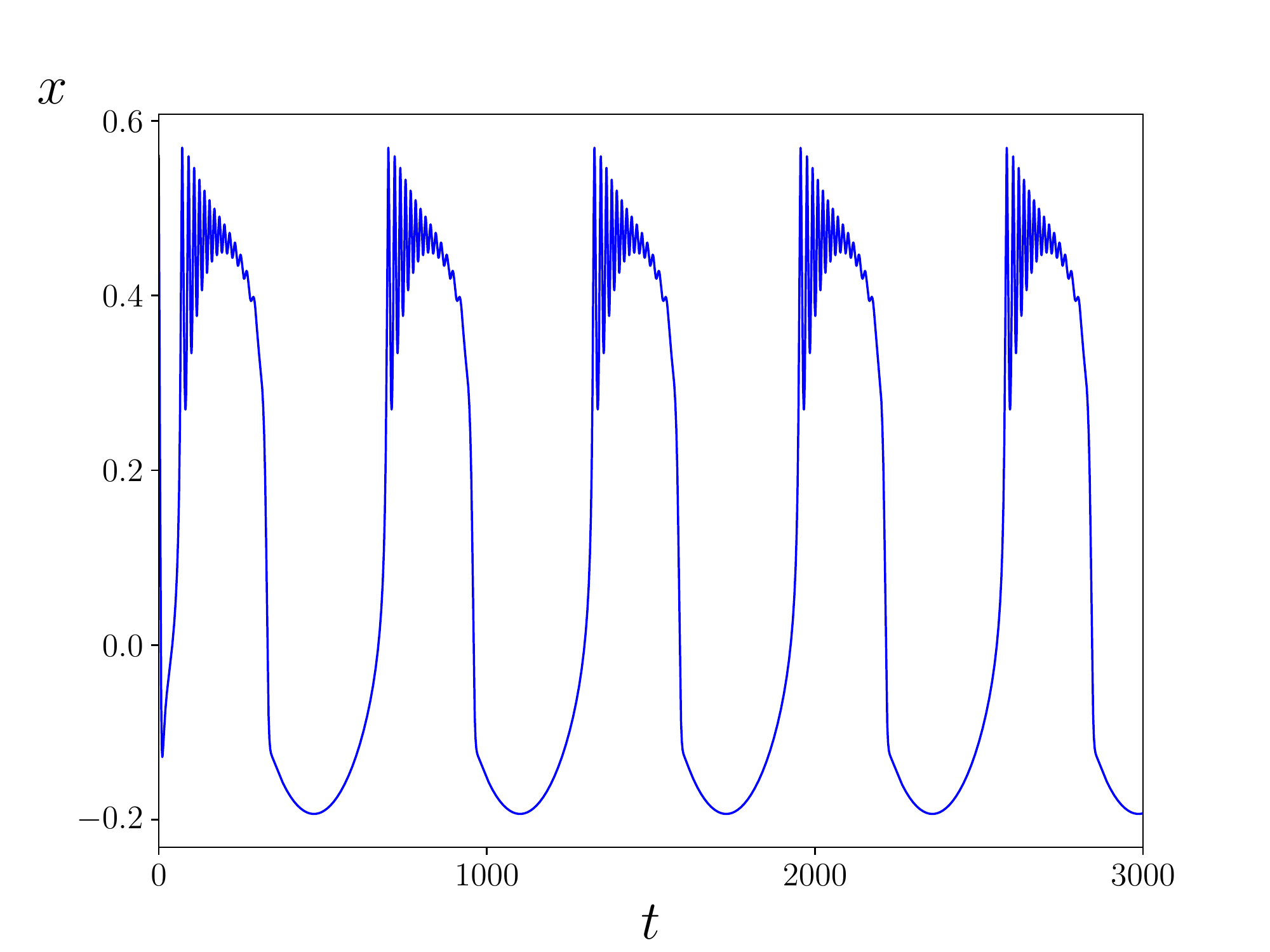}
  \label{fig:0.04}
\end{subfigure}%
\begin{subfigure}[b]{.26\textwidth}
  \centering
  \caption{}
 \includegraphics[width = \textwidth]{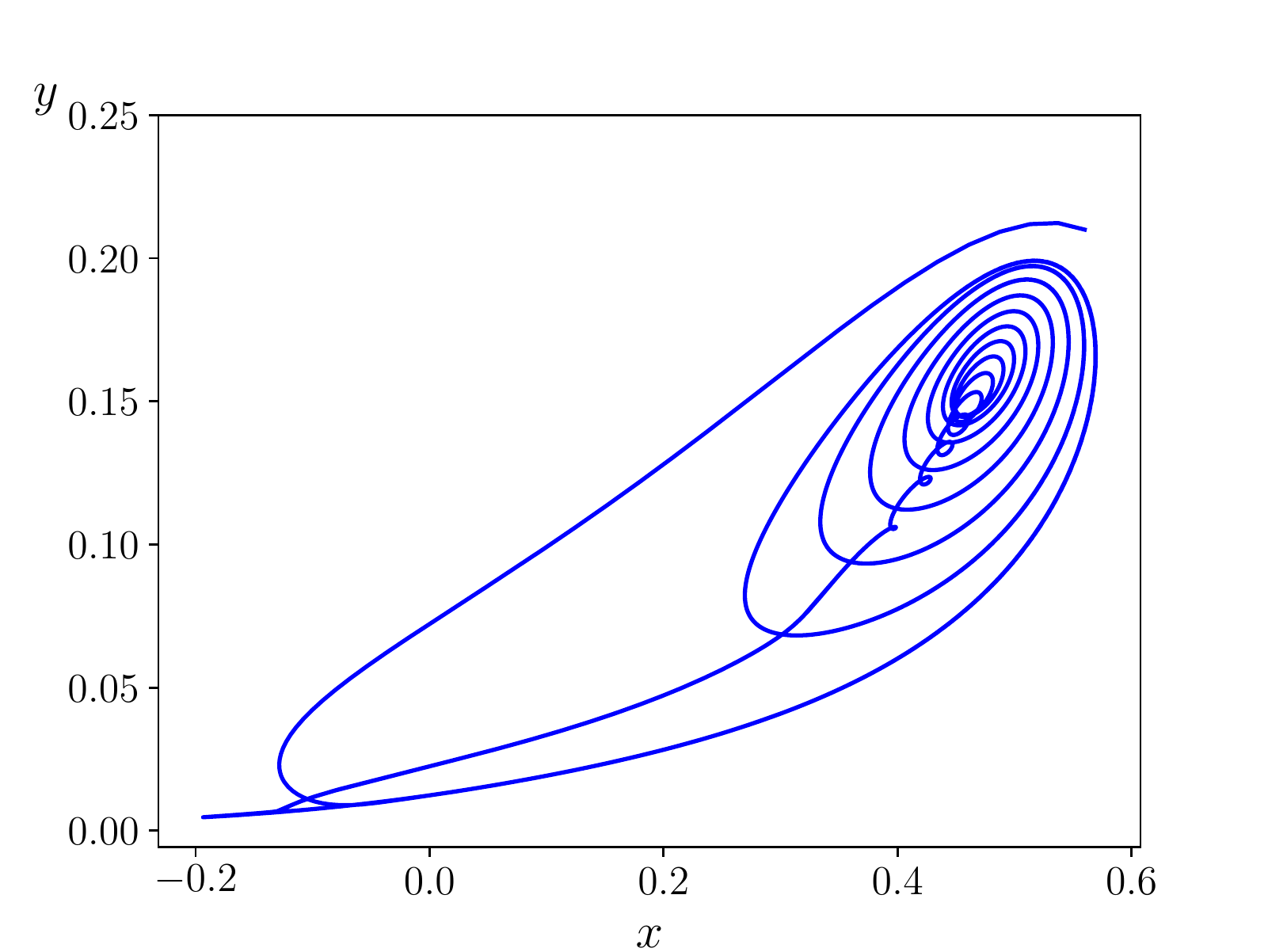}
  \label{fig:0.04_pp}
\end{subfigure}%
\caption{Time series plots of $x$ for model \eqref{eq:DML_em} for some values of $I_{0}$ with $\gamma=0.315$. (a) $I_{0}=0.00072$; (c) $I_{0}=0.0155$; (e) $I_{0}=0.016$; (g) $I_{0}=0.02$; (i) $I_{0}=0.04$. Their corresponding phase plane in $(x,y)$-plane are shown in (b), (d), (f), (h), and (j), respectively.}
\label{fig:TS_2}
\end{figure}

In this section, the dynamics of the improved denatured ML model \eqref{eq:DML_em} under the influence of electromagnetic induction and  external periodic current is investigated . The numerical simulations were carried out in {\sc Python} using odeint solver, which uses a variable step-size Runge-Kutta method. The time evolution of the membrane potential $x$ and the projection of the phase space in $(x,y)$-plane are shown in Fig.~\ref{fig:TS_2}. Various firing activities including regular spiking and bursting are found as the amplitude of the external current $I_{0}$ is varied and set $\gamma=0.315$. For $I_{0}=0.00072$, the model is in a rest state (see Fig.~\ref{fig:0.00072}). With the increase of $I_{0}$, regular periodic spiking activity is observed, an example is shown in Fig.~\ref{fig:0.0155} for $I_{0}=0.0155$. However, increasing $I_{0}$ further the model shows regular periodic bursting. The model has period-2, period-4 and period-13 when $I_{0}=0.016, 0.02,$ and $0.04$, respectively (see Figs.~\ref{fig:TS_2}e, g and i). As seen from the above analysis, complex dynamical behaviours occur in the improved denatured ML model under the influence of the amplitude of the external current $I_{0}$.

\begin{figure}[htbp]
\centering
\begin{subfigure}[b]{.26\textwidth}
  \centering
  \caption{}
  \includegraphics[width =\textwidth]{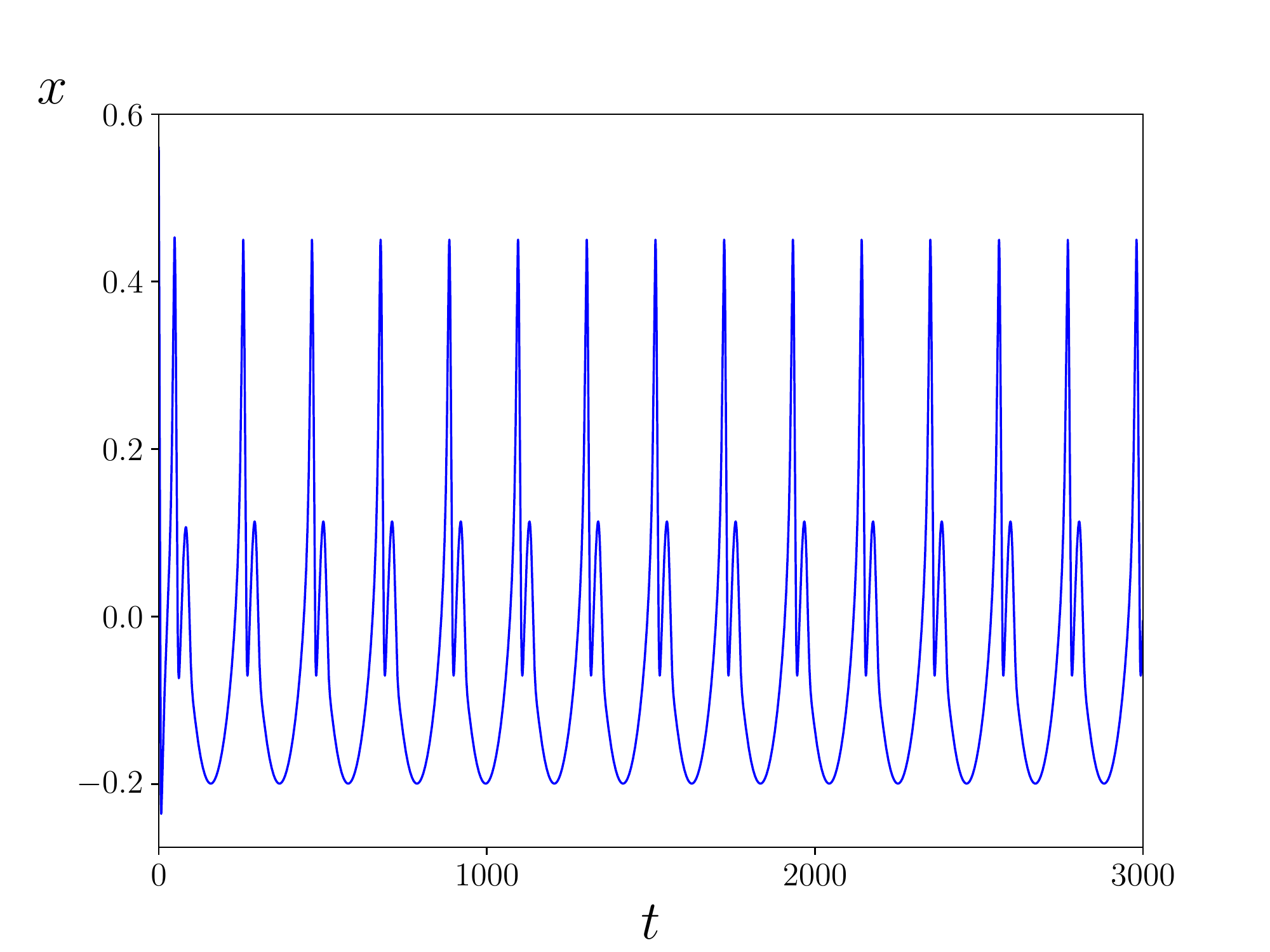}
\end{subfigure}%
\begin{subfigure}[b]{.26\textwidth}
  \centering
  \caption{}
  \includegraphics[width =\textwidth]{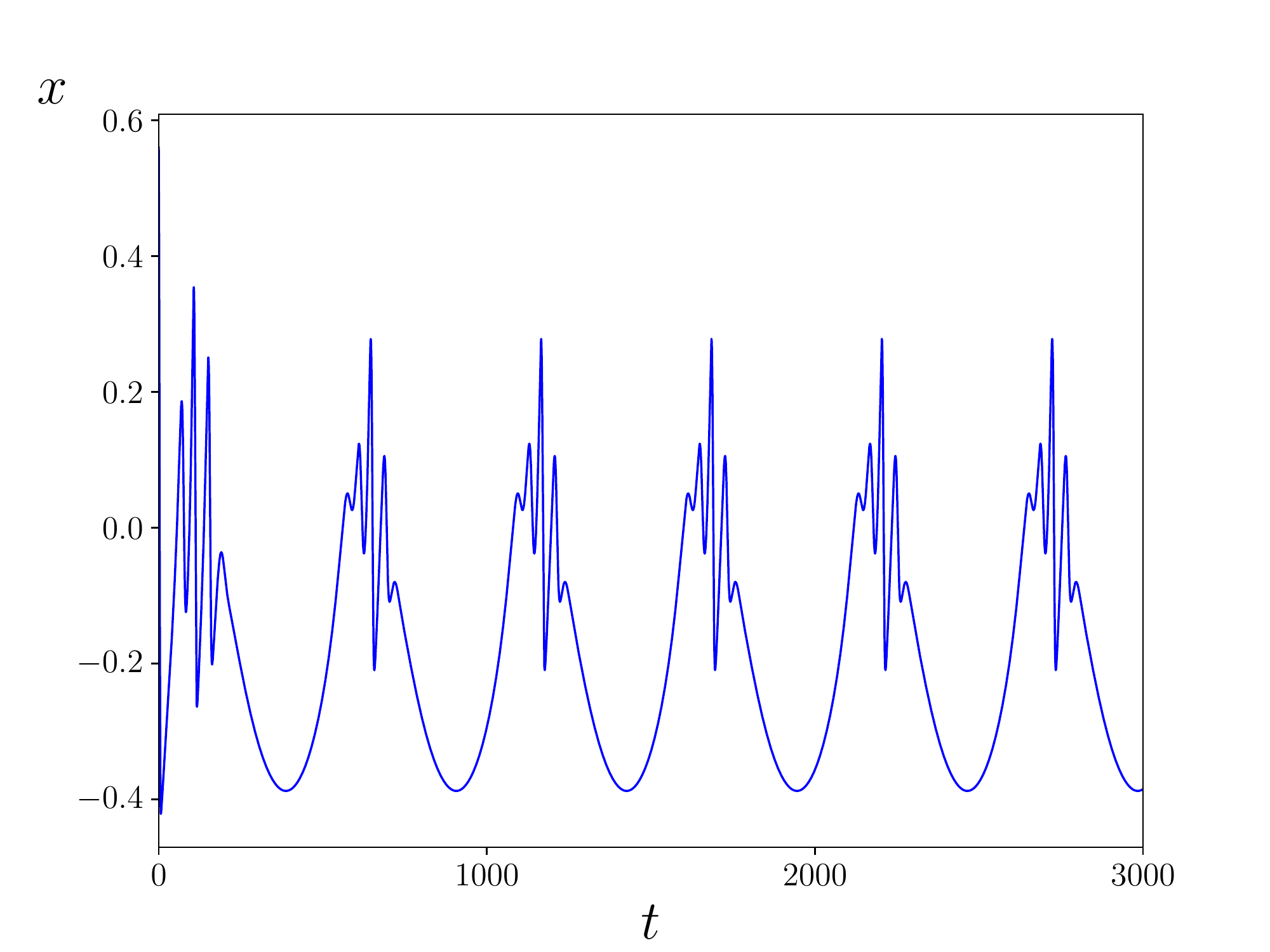}
\end{subfigure}%
\begin{subfigure}[b]{.26\linewidth}
    \centering
    \caption{}
    \includegraphics[width=\textwidth]{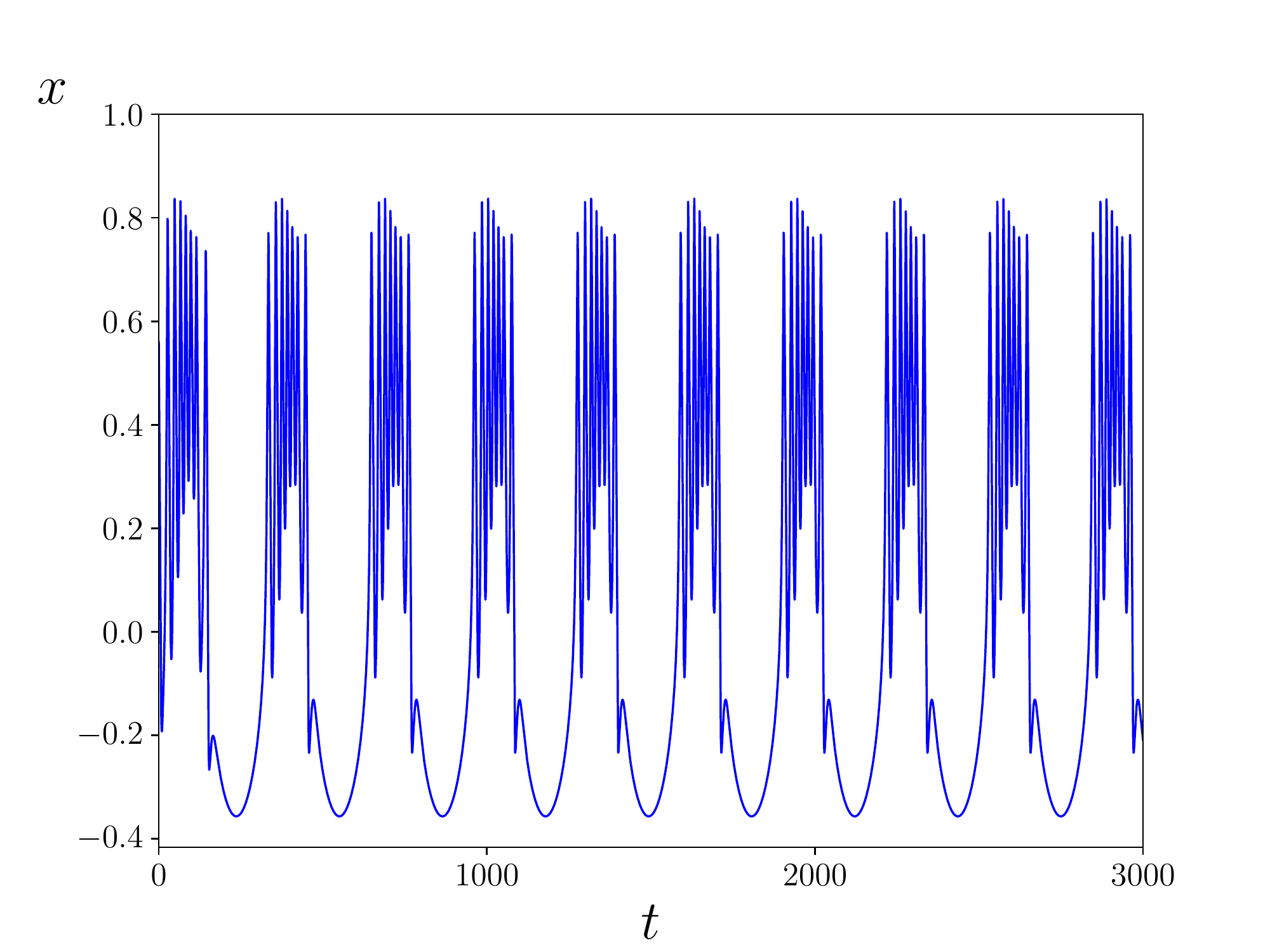}
  \end{subfigure}\\%
  \begin{subfigure}[b]{.26\linewidth}
    \centering
    \caption{}
    \includegraphics[width=\textwidth]{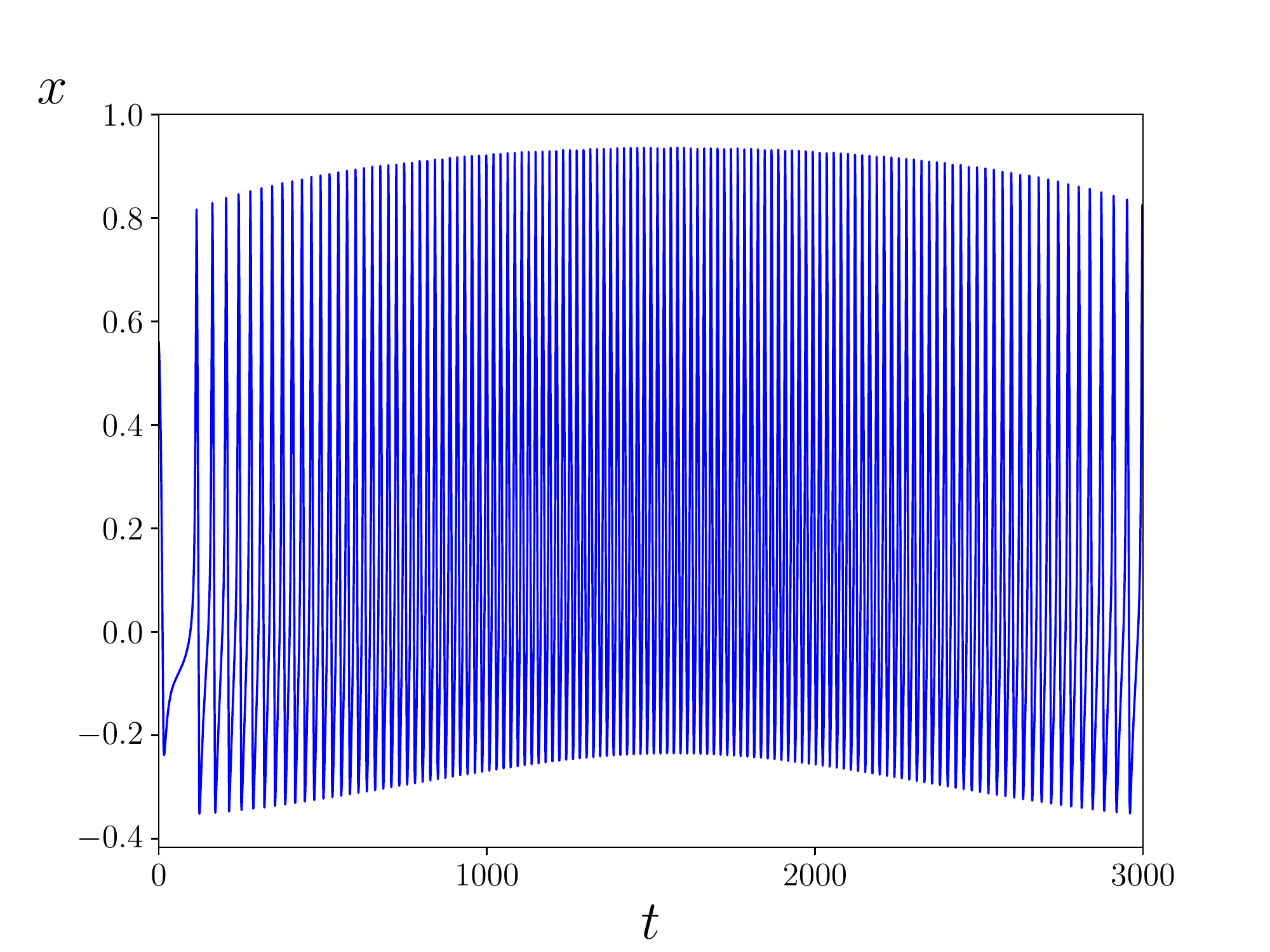}
  \end{subfigure}%
  \begin{subfigure}[b]{.26\textwidth}
  \centering
  \caption{}
 \includegraphics[width =\textwidth]{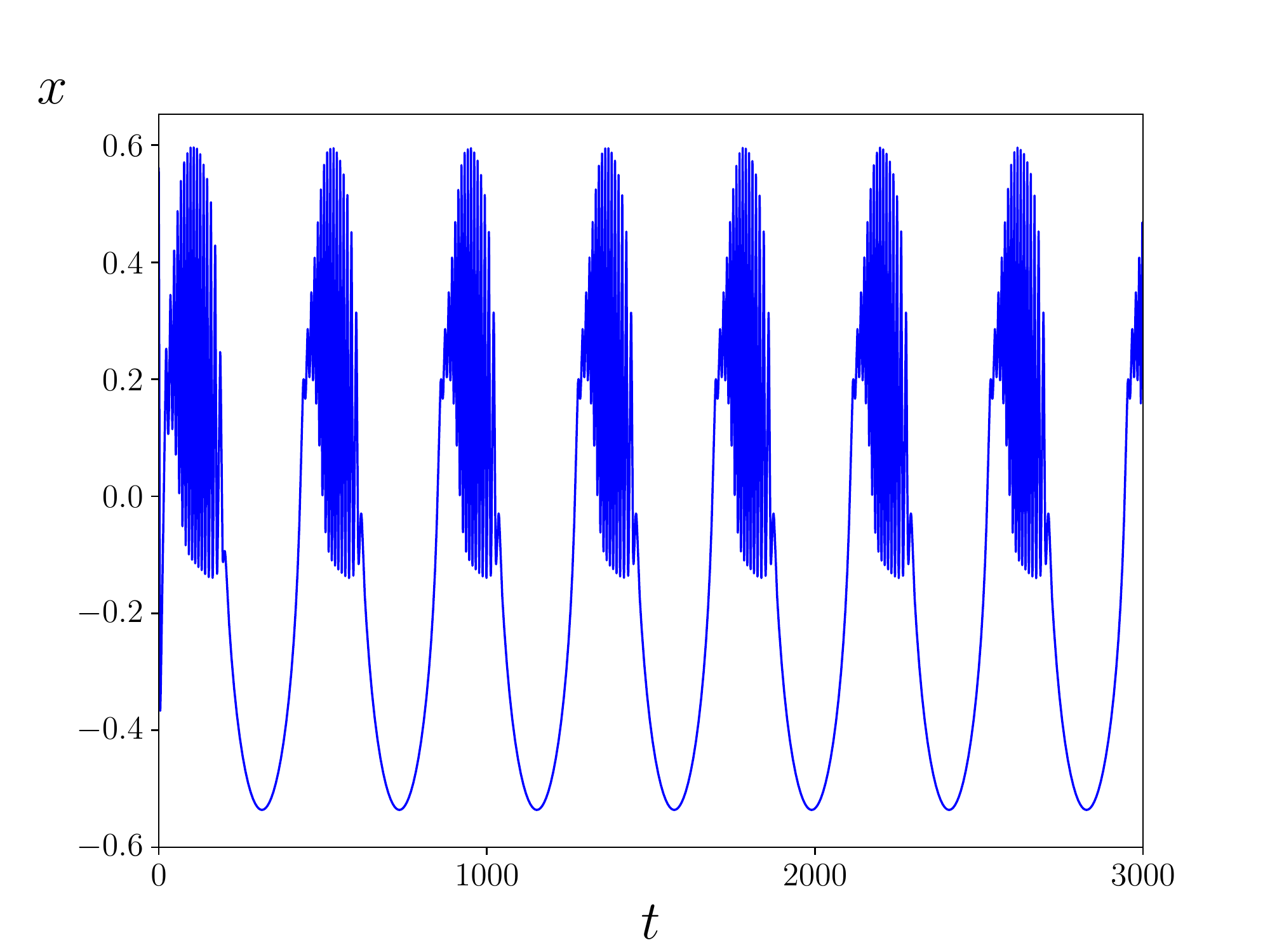}
\end{subfigure}%
\begin{subfigure}[b]{.26\textwidth}
  \centering
  \caption{}
 \includegraphics[width = \textwidth]{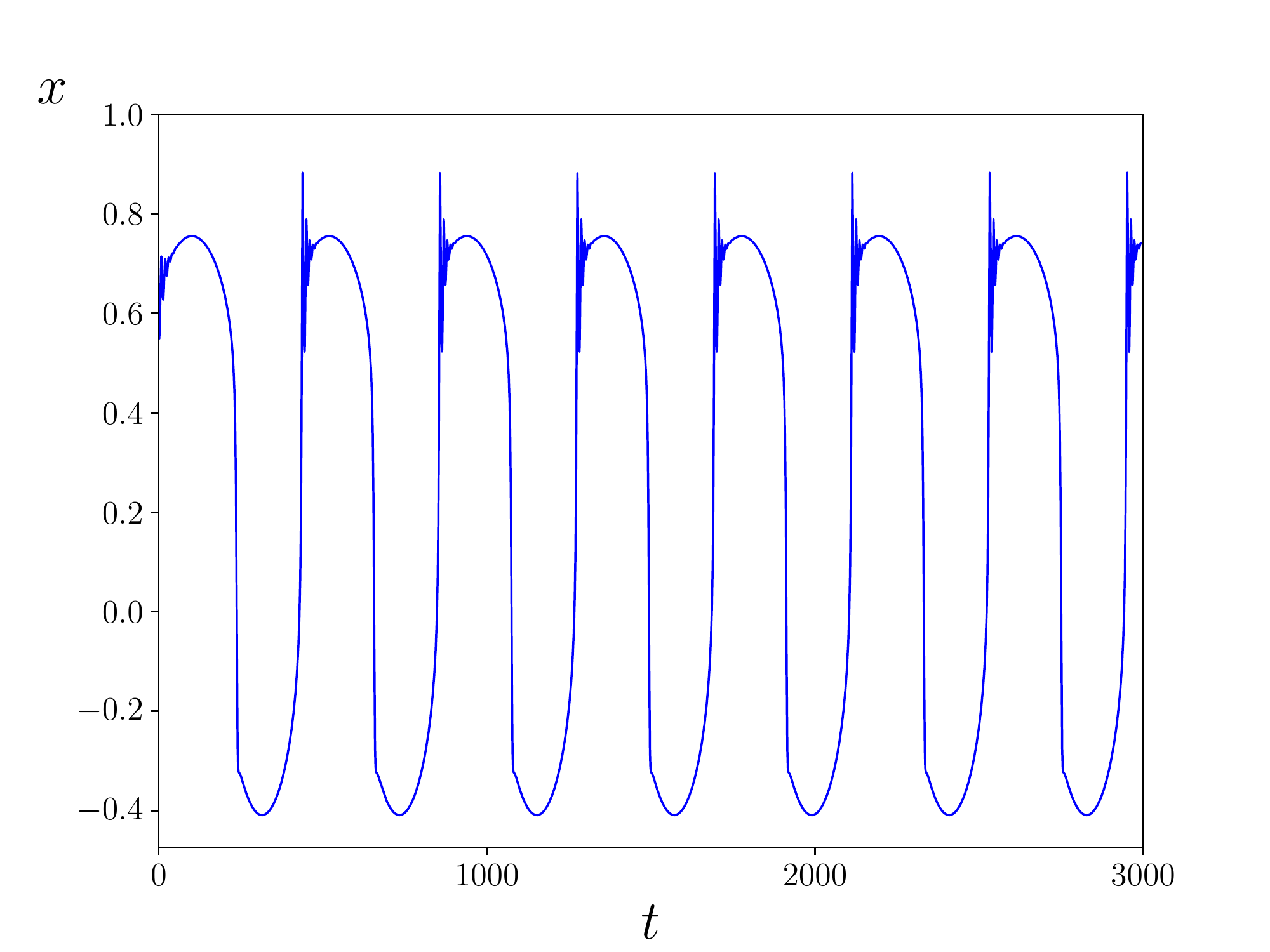}
\end{subfigure}\\%
\caption{Time series plots of $x$ for model \eqref{eq:DML_em} for (a) $A = 0.005, \gamma = 0.23, \omega = 0.03$, and $I_0 = 0.04$; (b) $A = 0.005, \gamma = 0.035, \omega = 0.0121, k = 0.0231$, and $I_0 = 0.187$; (c) $A = 0.002, \gamma = 0.1576, \omega = 0.02, k = 0.018$, and $I_0 = 0.17$; (d) $A = 0.002, \gamma = 0.1576, \omega = 0.001, k = 0.919$, and $I_0 = 0.17$; (e) $A = 0.0187, \gamma = 0.231, \omega = 0.015, k = 0.05$, and $I_0 = 0.434$; (f) $A = 0.0018, \gamma = 0.231, \omega = 0.015, k = 0.08$, and $I_0 = 0.201$.}
\label{fig:TS_3}
\end{figure}

Next, multiple parameters are varied simultaneously to explore the effects of the electromagnetic radiation on the dynamical behaviors of the improved denatured ML model. As seen in Fig.~\ref{fig:TS_3}, model \eqref{eq:DML_em} exhibits mixed mode oscillations and bursting activities. In Fig.~\ref{fig:TS_3}a-b, we observe mixed mode oscillations where several regular single spikes and subthreshold oscillations coexist. Both of these cases have $A=0.005$ with variations in $\gamma, \omega, k$ and $I_0$. For $A = 0.002, \gamma = 0.1576, \omega = 0.02, k = 0.018$, and $I_0 = 0.17$, there exists period-7 bursting with subthreshold oscillations (see Fig.~\ref{fig:TS_3}c). In Fig.~\ref{fig:TS_3}d, we observe that when $k = 0.919$, $\omega = 0.001$ with $A=0.002$, $\gamma = 0.1576$ , and $I_0 = 0.17$ the dynamics subsides to periodic spiking. Finally in Fig.~\ref{fig:TS_3}e-f regular periodic bursting patterns appear. 

\section{Conclusion}\label{sec:conclusion}
In this study, phase plane analysis and bifurcation analysis are employed to study the firing activities in a simple electrophysiological model of neurons. By using codimension-1 and -2 bifurcation analysis, we showed how variation of model parameters affects the dynamics of the membrane potential such as transitions from a rest state to periodic oscillations and vice versa. We found that model exhibits class II excitability. The class of excitability affects how neurons decode and process information, thus it is important to understand the mechanisms underlying the transitions between classes of excitability. As in \cite{Schaeffer}, the denatured ML model exhibits class I excitability. In future work we will investigate the mechanisms that govern switches between classes of excitability in the improved model.

In addition, we investigated the dynamics of the improved model behaviour under the influence of external periodic current and electromagnetic induction. Evidence from numerical simulations shows that the electrical activity of the improved neuron model depends greatly on system parameters. We observed that for different values of the external periodic current the model produces simple and complex oscillations including multiple mode bursting activities. Also, we explored the effects of the angular frequency $\omega$ and the feedback gain $k$ on the behaviour of the improved model.

The model displays qualitative behaviours observed in experiments and computational analysis of higher dimension neuron models. This work provide for some applications where simple models are required for physiological and pathophysiological responses in  neurons and other excitable cells.

In the present study, we have considered analysis of a single neuron, however, neurons are interconnected through chemical and electrical synapes. Hence, the analysis of ring network of neurons under electromagnetic induction will be a subject of future work.

  



\end{document}